\documentclass[12pt]{article}
\usepackage{amsfonts}

\usepackage{latexsym}
\usepackage{amssymb}
\usepackage{epsfig}
\usepackage{amsmath}
\usepackage{amsthm}
\usepackage[mathscr]{eucal}
\usepackage{subfigure}
\usepackage{graphicx}

\newtheorem{Theorem}{Theorem}

\newtheorem{Corollary}{Corollary}

\renewcommand{\qed}{\hfill{\ \ \rule{2mm}{2mm}} \vspace{0.2in}}

\newcommand{\ind}{1\hspace{-2.3mm}{1}}

\begin{document}

\title{Recurrence region of multiuser Aloha}
\author{ \textbf{Ghurumuruhan Ganesan}
\thanks{E-Mail: \texttt{gganesan82@gmail.com} } \\
\ \\
New York University, Abu Dhabi}
\date{}
\maketitle

\begin{abstract}
In this paper, we provide upper and lower bounds for the region of positive recurrence for a general finite user Aloha network.

\vspace{0.1in} \noindent \textbf{Key words:} Recurrence region, multiuser Aloha.

\vspace{0.1in} \noindent \textbf{AMS 2000 Subject Classification:} Primary:
60J10, 60K35; Secondary: 60C05, 62E10, 90B15, 91D30.
\end{abstract}

\bigskip

\setcounter{equation}{0}
\renewcommand\theequation{\thesection.\arabic{equation}}
\section{Introduction} \label{intro}
Consider a network of~\(M\) users assigned to access a common channel via random access. Much of previous literature has focused on finding bounds on the stability region for various channel models (see for example~\cite{naware,szpan,luo} and references therein). It is also important to determine the region where the overall network is positive recurrent. In this paper, we find bounds for the recurrence region for a general multiuser Aloha network. 

\subsection*{Channel Model}
The users are labelled~\(U_1,\ldots,U_M.\) Time is divided into slots of unit length and in each time slot, packets arrive randomly at the queue of each user. For~\(1 \leq i \leq M,\) let~\(A_i(n)\) be the random number of packets arriving at user~\(U_i\) in time slot~\(n.\) We assume that~\(\{A_i(n)\}_{n \geq 1}\) are independent and identically distributed (i.i.d) with
\begin{equation}\label{a_def}
\lambda_i = \mathbb{E}A_i(n)
\end{equation}
and satisfying
\begin{equation}\label{a1_cond}
\mathbb{P}(A_i(n) = 1) > 0.
\end{equation}
For~\(i \neq j,\) we also assume that~\(\{A_i(n)\}_{n \geq 1}\) is independent of~\(\{A_j(n)\}_{n \geq 1}.\) Depending on the queue length and the instantaneous channel conditions, user~\(U_i\) attempts to transmit packets. Formally, let~\(\{W_i(n)\}_{n \geq 1}\) be i.i.d nonnegative integer valued random variables with mean
\begin{equation}\label{w_def}
C_i = \mathbb{E}W_i(n)
\end{equation}
and
\begin{equation}\label{w1_cond}
\mathbb{P}(W_i(n) = 1) > 0.
\end{equation}

Let~\(Q_i(n)\) be the queue length of user~\(U_i\) at time slot~\(n.\) For time slot~\(n \geq 0,\) the update equation for~\(Q_i(.)\) is
\begin{eqnarray}
Q_i(n+1) &=& Q_i(n) + A_i(n+1) \nonumber\\
&&\;\;\;\;\;-\;\;\min(Q_i(n),W_i(n+1))\ind\left(B_i(n+1) \bigcap \bigcap_{j \neq i} B^c_j(n+1)\right) \nonumber\\
\label{q_def}
\end{eqnarray}
where
\begin{equation}\label{b_def}
B_i(n+1) = \left\{Q_i(n) \geq 1\right\} \bigcap \left\{W_i(n+1) \geq 1\right\}
\end{equation}
is the intersection of the events that the queue of user~\(U_i\) is nonempty and the user~\(U_i\) attempts to transmit in time slot~\(n+1.\) Define
\begin{equation}\label{p_def}
p_i := \mathbb{P}(W_i(n) \geq 1)
\end{equation}
to be the attempt probability of user~\(U_i.\)

Conditions~(\ref{a1_cond}) and~(\ref{w1_cond}) ensure that the Markov chain \[\underline{Q}(n) = (Q_1(n),\ldots,Q_M(n)), n \geq 1\] is irreducible. The following result provides necessary and sufficient conditions for recurrence of~\(\{\underline{Q}(n)\}.\)
\begin{Theorem}\label{thm_main}
If~\(W_i(n) \in \{0,1\}\) for all~\(n\) and
\begin{equation}\label{eq_low}
\sum_{i=1}^{M} \frac{\lambda_i}{C_i\prod_{j\neq i}(1-p_j)} < 1,
\end{equation}
then~\(\{\underline{Q}(n)\}\) is positive recurrent.

If~\(\max_{1 \leq i \leq M} \mathbb{E}A_i^4(1)  < \infty, \max_{1 \leq i \leq M} \mathbb{E}W_i^4(1) < \infty\) and
\begin{equation}\label{eq_up}
\lambda_i  > C_i \prod_{j \neq i} (1-p_j)
\end{equation}
for all~\(1 \leq i \leq M,\) then~\(\{\underline{Q}(n)\}\) is transient.
\end{Theorem}

For the particular case when~\(W_i(n) \in \{0,1\},\) we obtain the usual single packet Aloha network with
\begin{equation}\label{p_def2}
p_i := \mathbb{P}(W_i(n) =1) = 1 - \mathbb{P}(W_i(n) = 0).
\end{equation}
For~\(\underline{p} = (p_1,\ldots,p_M),\) the Markov chain~\(\underline{Q}(n) = \underline{Q}(\underline{p},n)\)  and so define the recurrence region
\begin{eqnarray}
{\cal R} &=& \left\{(\lambda_1,\ldots,\lambda_M) : \exists\;\underline{p} \in (0,1)^{M} \text{ such that } \right. \nonumber\\
&&\;\;\;\; \left. \{\underline{Q}(\underline{p},n)\}_{n \geq 1} \text{ is positive recurrent}\right\}\label{stab_reg_def}
\end{eqnarray}

To obtain bounds for~\({\cal R}\) using Theorem~\ref{thm_main}, we have some definitions. For~\(\underline{p} = (p_1,\ldots,p_M) \in (0,1)^{M},\) define
\begin{equation}\label{st_def11}
{\cal C}_1(\underline{p}) := \left\{(\lambda_1,\ldots,\lambda_M) \in (0,1)^{M} : \sum_{i=1}^{M} \frac{\lambda_i}{p_i \prod_{j \neq i} (1-p_j)} < 1\right\}
\end{equation}
and let
\begin{equation}\label{st_def_low}
{\cal C}_1 = \bigcup_{ \underline{p} \in (0,1)^{M} } {\cal C}_1(\underline{p}).
\end{equation}
Similarly define
\begin{equation}\label{st_def1}
{\cal C}_2(\underline{p}) := \left\{(\lambda_1,\ldots,\lambda_M) \in (0,1)^{M} : \lambda_i < p_i \prod_{j \neq i}(1-p_j) \text{ for some }1 \leq i \leq M\right\}
\end{equation}
and let
\begin{equation}\label{st_def2}
{\cal C}_2 = \bigcup_{\underline{p} \in (0,1)^{M}} {\cal C}_2(\underline{p}).
\end{equation}
\begin{Corollary} \label{thm1} The recurrence region~\({\cal R}\) defined in~(\ref{stab_reg_def}) satisfies
\begin{equation}\label{st_main}
{\cal C}_1 \subseteq {\cal R} \subseteq \text{cl}({\cal C}_2)
\end{equation}
where~\(\text{cl}({\cal C}_2)\) is the closure of the set~\({\cal C}_2.\)
\end{Corollary}




The paper is organized as follows: In Section~\ref{pf1}, we prove Theorem~\ref{thm1}.

\section*{Proof of Theorem~\ref{thm1}}\label{pf1}
\subsection*{Lower bound}
From~(\ref{eq_low}), there exists~\(\epsilon_0 > 0\) so that
\begin{equation}\label{eq_lam}
\sum_{i=1}^{M} \frac{\lambda_i}{v_i} \leq 1-\epsilon_0,
\end{equation}
where~\(v_i := C_i \prod_{j \neq i}(1-p_j).\)

For~\(1 \leq i \leq M\) we have from the definition of the event~\(B_i(.)\) in~(\ref{b_def}) that \[\ind(B_j(n+1)) \leq \ind(W_j(n+1) \geq 1)\] for all~\( j \neq i.\) Moreover if~\(Q_i(n) \geq 1,\) then~\(\min(Q_i(n),W_i(n+1)) = W_i(n+1)\) since~\(W_i(n+1) \in \{0,1\}\) and so from~(\ref{q_def}),
\begin{equation}\label{q1}
Q_i(n+1) \leq Q_i(n)  + A_i(n+1) - W_i(n+1) \ind(Q_i(n) \geq 1)\ind(V_i(n+1)),
\end{equation}
where
\begin{equation}\label{v_def}
V_i(n+1) := \{W_i(n+1) \geq 1\} \bigcap \bigcap_{j \neq i} \{W_j(n+1) = 0\}
\end{equation}
is the event that only user~\(U_i\) attempts to transmit at time slot~\(n+1.\) The term~\(W_i(n+1) \ind(V_i(n+1))\) is independent of~\(Q_i(n)\) and has mean
\begin{equation}\label{vi_exp}
\mathbb{E}W_i(n+1)\ind(V_i(n+1)) = C_i \prod_{j \neq i}(1-p_j) = v_i,
\end{equation}
from~(\ref{eq_lam}).

Let~\(Q_i(0) = 0\) for all~\(1 \leq i \leq M\) and define
\begin{equation}\label{t_def}
T = \inf\{k \geq 1 : Q_i(k)  = 0 \text{ for all }1 \leq i \leq M\}
\end{equation}
to be the first time that the queues of all the users are simultaneously empty, again. Analogous to~\cite{asmu}, we show that the expected time to return to the origin~\(\mathbb{E}T\) is finite. Multiplying both sides of~(\ref{q1}) by~\(\ind(T \geq n+1)\) and using the fact that~\(Q_i(n)\ind(T \geq n+1) \leq Q_i(n) \ind(T \geq n),\) we get
\begin{eqnarray}
&&Q_i(n+1)\ind(T \geq n+1) \nonumber\\
&&\;\;\;\leq Q_i(n)\ind(T \geq n) + A_i(n+1)\ind(T \geq n+1) \nonumber\\
&&\;\;\;\;\;\;-\;\;\ind(Q_i(n) \geq 1) \ind(T \geq n+1) W_i(n+1) \ind(V_i(n+1)). \label{q2}
\end{eqnarray}
Defining~
\begin{equation}\label{y_n_i_def}
y_n(i) := \mathbb{E}Q_i(n)\ind(T \geq n)
\end{equation}
we get from~(\ref{q2}) that
\begin{eqnarray}
y_{n+1}(i) &\leq& y_{n}(i) + \mathbb{E}A_i(n+1)\ind(T \geq n+1)  \nonumber\\
&&\;\;\;\;\;\;-\;\;\mathbb{E}\ind(Q_i(n) \geq 1)\ind(T \geq n+1) W_i(n+1)\ind(V_i(n+1)) \nonumber\\
&=& y_{n}(i) + \mathbb{E}A_i(n+1)\mathbb{E}\ind(T \geq n+1) \nonumber\\
&&\;\;\;\;\;\;-\;\;\mathbb{E}\ind(Q_i(n) \geq 1) \ind(T \geq n+1) \mathbb{E}W_i(n+1)\ind(V_i(n+1)) \nonumber\\
\label{gt1}\\
&=& y_n(i) + \lambda_i \mathbb{P}(T \geq n+1) - v_i \mathbb{P}(\{Q_i(n) \geq 1\} \cap \{T \geq n+1\}). \nonumber\\
\label{gt}
\end{eqnarray}
The relation~(\ref{gt1}) follows from the fact~\(\{T \geq n+1\} = \{T \leq n\}^{c}\) and that~\(A_i(n+1), W_i(n+1)\) and~\(V_i(n+1)\) are independent of the process up to~\(n\) time slots. The final estimate in~(\ref{gt}) is obtained using~(\ref{vi_exp}).

Defining
\begin{equation}\label{yn_deff}
y_n := \sum_{i=1}^{M} \frac{y_n(i)}{v_i},
\end{equation}
we get from~(\ref{gt}) that
\begin{equation}\label{yn_rec}
y_{n+1} \leq y_n + \left(\sum_{i=1}^{M}\frac{\lambda_i}{v_i}\right)\mathbb{P}(T \geq n+1) - \Delta_n
\end{equation}
where
\begin{eqnarray}
\Delta_n &:=& \sum_{i=1}^{M}\mathbb{P}\left(\{Q_i(n) \geq 1\} \bigcap \{T \geq n+1\}\right) \nonumber\\
&\geq& \mathbb{P}\left(\bigcup_{i=1}^{M} \{Q_i(n) \geq 1\} \bigcap \{T \geq n+1\}\right) \nonumber\\
&=& \mathbb{P}\left(T \geq n+1\right), \label{del_n}
\end{eqnarray}
since if~\(T \geq n+1,\) then at least one of the queues at time slot~\(n \geq 1\) is non empty.

Using~(\ref{del_n}) in~(\ref{yn_rec}) gives
\begin{eqnarray}
y_{n+1} &\leq& y_n + \left(\sum_{i=1}^{M} \frac{\lambda_i}{v_i}-  1\right) \mathbb{P}(T \geq n+1) \nonumber\\
&\leq& y_n - \epsilon_0 \mathbb{P}(T \geq n+1) \nonumber
\end{eqnarray}
by~(\ref{eq_lam}). Thus
\begin{equation}\nonumber
\mathbb{P}(T \geq n+1) \leq \frac{1}{\epsilon_0}(y_n - y_{n+1})
\end{equation}
and adding telescopically gives for~\(J \geq 1\) that
\begin{equation}\label{t_fin}
\sum_{k=1}^{J}\mathbb{P}(T \geq k+1) \leq \frac{1}{\epsilon_0}(y_1 - y_{J+1}) \leq \frac{1}{\epsilon_0}y_1,
\end{equation}
where~\[y_1 = \sum_{i=1}^{M} \frac{y_1(i)}{v_i}\] using~(\ref{yn_deff}) and~\[y_1(i) = \mathbb{E}Q_i(1) \ind(T \geq 1) \leq \mathbb{E}Q_i(1) \leq \mathbb{E}A_i(1) = \lambda_i < \infty\] using~(\ref{y_n_i_def}). This implies that~\(0 \leq y_1 < \infty\) and since~\(J\) is arbitrary, we get from~(\ref{t_fin}) that~\(\mathbb{E}T < \infty\) and so the Markov chain~\(\{\underline{Q}(n)\}\) is positive recurrent.



\subsection*{Upper bound}
Let~\(K \geq 1\) be a large integer constant to be determined later and let~\(T_0 := 0\) and~\(T_1 := K.\)
We now observe the overall queue process from time slot\\\(T_0+1 = 1\) to time slot~\(T_1.\)  Recall from~(\ref{q_def}) that~\(W_i(n+1)\) is the maximum number of packets transmitted by user~\(U_i\) in time slot~\(n+1.\) If
\begin{equation}\label{z_i_def}
Z_i(T_0,T_1) := \left\{\sum_{n = T_0}^{T_1-1} W_i(n+1) < 2C_i(T_1-T_0) \right\}
\end{equation}
then
\begin{equation}\label{z_it_est}
\mathbb{P}\left(Z^c_i(T_0,T_1)\right) = \mathbb{P}\left(S_i \geq C_i(T_1-T_0)\right),
\end{equation}
where~\(S_i = \sum_{n=T_0}^{T_1-1} (W_i(n+1) - C_i) \) is a sum of independent zero mean random variables and so
\begin{equation}\label{es42}
\mathbb{E}S^4_i = \sum_{n} \mathbb{E}(W_i(n+1)-C_i)^4 + \sum_{n \neq m} \mathbb{E}(W_i(n+1)-C_i)^2\mathbb{E}(W_i(m+1)-C_i)^2.
\end{equation}

Using the finite fourth moment condition of~\(W_i(n)\) (see statement prior to~(\ref{eq_up})), the first term in~(\ref{es42}) is~\(\mathbb{E}(W_i(1) - C_i)^4 (T_1-T_0)\) and the second term in~(\ref{es42}) is at most~\[(T_1-T_0)^2(\mathbb{E}(W_i(1)-C_i)^2)^2 \leq (T_1-T_0)^2\mathbb{E}(W_i(1)-C_i)^4.\] Combining,
\begin{equation}
\mathbb{E}S^4_i \leq \alpha_1(i) (T_1-T_0)^2 \label{es4}
\end{equation}
for some constant~\(\alpha_1(i) > 0,\) not depending on~\(T_1\) or~\(T_0.\) From~(\ref{z_it_est}),~(\ref{es4}) and Markov inequality, we get
\begin{equation}\label{zit2}
\mathbb{P}\left(Z^c_i(T_0,T_1)\right) \leq \frac{\alpha_1(i)(T_1-T_0)^2}{C_i^4(T_1-T_0)^4} \leq \frac{\alpha_2}{(T_1-T_0)^2}
\end{equation}
where~\(\alpha_2 = \max_{1 \leq i \leq M} \frac{\alpha_1(i)}{C_i^4}\) is a constant. If
\begin{equation}\label{z_t_def}
Z(T_0,T_1) := \bigcap_{i=1}^{M} Z_i(T_0,T_1),
\end{equation}
then from~(\ref{zit2}),
\begin{equation}\label{zt_est2}
\mathbb{P}\left(Z(T_0,T_1)\right) \geq 1-\frac{\alpha_2 M}{(T_1-T_0)^2}.
\end{equation}

For~\(1 \leq i \leq M,\) set the initial queue length
\begin{equation}\label{q_i_0}
Q_i(T_0) = Q_i(0) = \delta (T_1-T_0) =: L_1,
\end{equation}
where~\(\delta = 3\max_{1 \leq i \leq M} C_i\) and~\(T_1 = K\) is large so that~\(\delta(T_1 -T_0) = \delta K  > 1.\) If~\(Z(T_0,T_1)\) occurs, then at most~\(2C_i(T_1-T_0)\) packets are transmitted from user~\(U_i\) and so the queue of user~\(U_i\) never becomes empty between time slots~\(T_0+1\) and~\(T_1.\) From the queue update equation~(\ref{q_def}) we therefore get for~\(T_0 \leq n \leq T_1\) and~\(1 \leq i \leq M\) that
\begin{equation}\label{qw_def}
Q_i(n+1) \geq Q_i(n)  + A_i(n+1) - W_i(n+1)\ind(V_i(n+1))
\end{equation}
where~\(V_i(.)\) is as defined in~(\ref{v_def}). Adding telescopically,
\begin{equation}\label{qw_teles}
Q_i(T_1) = Q_i(T_0)  + R_i(T_0,T_1),
\end{equation}
where
\begin{equation}\label{rk_def}
R_i(T_0,T_1) := \sum_{n=T_0}^{T_1-1}(A_i(n+1) - W_i(n+1)\ind(V_i(n+1))).
\end{equation}
From~(\ref{eq_up}) we have that
\begin{equation}\label{re_est}
\mathbb{E}(A_i(n+1) - W_i(n+1)\ind(V_i(n+1))) = 2\epsilon_1(i) > 0
\end{equation}
for all~\(n.\) Moreover, the term~\(R_i(T_0,T_1)\) is also a sum of i.i.d zero mean random variables and so arguing as in~(\ref{zit2}), we get
\begin{equation}\label{rk_est2}
\mathbb{P}\left(R_i(T_0,T_1) \geq \epsilon_1(i) (T_1-T_0) \right) \geq 1-\frac{\alpha_3(i)}{(T_1-T_0)^2}
\end{equation}
for some constant~\(\alpha_3(i) > 0,\) not depending on~\(T_1\) or~\(T_0.\)

Letting~\(\epsilon_1 = \min_{1 \leq i \leq M}\epsilon_1(i)\) and~\(\alpha_3 = \max_{1 \leq i \leq M} \alpha_3(i)\) and defining
\begin{equation}\label{rt0_def}
X(T_0,T_1) := \bigcap_{i=1}^{M} \left\{R_i(T_0,T_1) \geq \epsilon_1(T_1-T_0)\right\},
\end{equation}
we get from~(\ref{rk_est2}) that
\begin{equation}\label{rk_est3}
\mathbb{P}\left(X(T_0,T_1) \right) \geq 1-\frac{\alpha_3 M}{(T_1-T_0)^2}
\end{equation}
and if
\begin{equation}\label{yt_def}
Y(T_0,T_1) := Z(T_0,T_1) \bigcap X(T_0,T_1),
\end{equation}
then  from~(\ref{z_it_est}) and~(\ref{rk_est3}), we get
\begin{equation}\label{y_est}
\mathbb{P}(Y(T_0,T_1)) \geq 1-\frac{\alpha_4}{(T_1-T_0)^2}
\end{equation}
for some constant~\(\alpha_4 > 0,\) not depending on~\(T_0\) or~\(T_1.\)

Suppose now that~\(Y(T_0,T_1)\) occurs. Between time slots~\(T_0\) and~\(T_1,\) none of the queues of the~\(M\) users ever becomes empty and at time slot~\(T_1,\) the queue length~\(Q_i(T_1)\) is at least~
\begin{equation}\label{l2_est}
Q_i(T_0) + \epsilon_1 (T_1-T_0) = L_1 + \epsilon_1 (T_1-T_0) = (\delta + \epsilon_1)(T_1-T_0) =: L_2,
\end{equation}
using~(\ref{q_i_0}) and~(\ref{qw_teles}). For~\(j \geq 2,\) we now repeat the above procedure between time slots~\(T_{j-1}+1\) and~\(T_j,\) where~\(T_j\) is determined by the relation
\begin{equation}\label{rec_lj}
\delta(T_j-T_{j-1}) = L_{j} = (\delta+\epsilon_1)(T_{j-1} - T_{j-2})
\end{equation}
and~\(\delta > 0\) is as in~(\ref{q_i_0}). Using the first and last relations in~(\ref{rec_lj}) iteratively, we get
\begin{equation}\label{lj_est}
L_j = \delta(T_1-T_0)\left(1 + \frac{\epsilon_1}{\delta}\right)^{j-1} \geq (T_1-T_0)(\delta + (j-1)\epsilon_1)
\end{equation}
and so from~(\ref{rec_lj}),
\begin{equation}\label{tj_diff}
T_j-T_{j-1} \geq (T_1-T_0)\left(1+ (j-1)\frac{\epsilon_1}{\delta}\right).
\end{equation}
Also analogous to~(\ref{y_est}), we have
\begin{equation}\label{y_estj}
\mathbb{P}(Y(T_{j-1},T_j)) \geq 1-\frac{\alpha_4}{(T_j-T_{j-1})^2} \geq 1-\frac{\alpha_5}{(\delta + (j-1)\epsilon_1)^2}
\end{equation}
for all~\(j \geq 2\) and for some constant~\(\alpha_5 > 0\) not depending on~\(j.\)

If the event
\begin{equation}\label{y_def_ov}
Y := \bigcap_{j \geq 1} Y(T_{j-1},T_j)
\end{equation}
occurs, then none of the queues of any user ever becomes empty. Using~(\ref{y_estj}) and the Markov property we also have
\[\mathbb{P}(Y) \geq \prod_{j} \left(1- \frac{\alpha_6}{(\delta + (j-1)\epsilon_1)^2}\right) > 0,\]
since~\(\sum_j \frac{\alpha_6}{(\delta + (j-1)\epsilon_1)^2} < \infty.\) Recall that the initial queue length of each user is~\(\delta(T_1-T_0) = \delta K\) (see~(\ref{q_i_0})) and so starting from~\((\delta K, \ldots,\delta K),\) the above discussion implies that with positive probability, the Markov chain~\(\{\underline{Q}(n)\}\) never reaches the origin. Since~\(\{\underline{Q}(n)\}\) is irreducible, this implies that starting from the origin, the chain~\(\{\underline{Q}(n)\}\) never returns to the origin, with positive probability. Therefore~\(\{\underline{Q}(n)\}\) is transient.~\(\qed\)








\bibliographystyle{plain}

\end{document}